\documentclass{amsart}
\usepackage[all]{xy}
\usepackage{graphicx}
\usepackage{psfrag}
\usepackage{amsthm}
\usepackage{amscd}
\usepackage{amsfonts}
\usepackage{amstext}
\usepackage{amssymb}
\usepackage{amsmath}
\usepackage{amsxtra}
\usepackage{plain}
\usepackage[all]{xy}

\usepackage[active]{srcltx}
\usepackage{colordvi}

\newtheorem{prop}{\bf Proposition}[section]
\newtheorem{cor}[prop]{{\bf Corollary}}
\newtheorem{lem}[prop]{{\bf Lemma}}
\newtheorem{thm}[prop]{{\bf Theorem}}

\numberwithin{equation}{section}

\newtheorem{rem}{{\bf Remark }}{ }

\newenvironment{pf}{{\bf Proof: }}{\qed\endtrivlist}

\newcommand{\Z}{\mathbb{Z} }
\newcommand{\N}{\mathbb{N} }
\newcommand{\C}{\mathbb{C} }
\newcommand{\R}{\mathbb{R} }

\newcommand{\tr}{\operatorname{tr}}

\newcommand{\ind}{\operatorname{ind}}

\newcommand{\T}{\mathbb{T}}

\begin{document}

\title{Analytic approach to $S^1$-equivariant Morse inequalities}
\author{Mostafa E. Zadeh  \and Reza Moghadasi}
\address{Mostafa Esfahani Zadeh, 
Sharif University of Technology,, Tehran-Iran
\newline
Reza Moghadasi
Sharif University of Technology,
Tehran-Iran}
\email{esfahani@sharif.ir}
\email{moghadasi@sharif.ir}

\begin{abstract}

The cohomology groups of a closed manifold $M$ can be reconstructed 
using the gradient flow of a Morse-Smale function $f\colon M\to \R$. 
A direct result of this 
construction are Morse inequalities that provide  
lower bounds for the number of critical points of $f$ in term of Betti numbers of $M$. 
E. Witten showed that these inequalities can be deduced analytically by studying the asymptotic behaviour  of the 
deformed Laplacian operator.  
In this paper, adopting Witten's approach, we provide an analytic proof for  
the so-called equivariant Morse inequalities when the underlying manifold is acted upon by 
the Lie group $G=S^1$ and the Morse function $f$ is invariant with respect to this action.  
\end{abstract}

\maketitle

\section{Introduction}

Classical Morse theory involves  the fundamental observation that the cellular 
structure of a closed manifold $M$ can be reconstructed through level sets of a 
Morse function $f$ on $M$. In particular this gives a way to 
reconstruct the cellular chain complex and therefore the cohomology of $M$, 
as is explained clearly in \cite{Milnor-Morse} and \cite{Schwarz-Morse}. 

The seminal paper of E. Witten \cite{Witten}, which was inspired 
by ideas from quantum field theory, shed a new light on Morse theory by providing 
a new chain complex for reconstructing the cohomology of $M$. 
This complex, called the Morse-Smale-Witten 
complex, is  generated by the critical points of $f$ and graded by their 
Morse indices, as in the cellular complex. However its  
differentials are defined using the gradient lines between critical points 
whose indices differ by one. 
We refer to \cite{Schwarz-Morse} for a detailed exposition of this theory. Amongst others, 
this construction led to the innovation of Floer Homology and solved (partially) 
the Arnold conjecture, c.f. \cite{Salamon-Lectures}. 

An immediate consequence of the reconstruction of the singular cohomology 
via critical points of a Morse function is the Morse inequalities. 
Roughly speaking, they provide lower bounds   
for the number of critical points in term of the Betti numbers of $M$, i.e. 
the rank of the cohomology groups of $M$. 
If one is interested in these inequalities rather than the cohomology itself, 
there is a very elegant and conceptual analytic derivation. This is  
Witten's idea of deforming the de Rham complex in an appropriate way using the 
Morse function and then studying the asymptotic behaviour of this complex. 
In this paper we follow Roe's account of Witten's approach in 
\cite[chapter 14]{Roe-elliptic}. 

Morse theory can be generalized in other directions. For instance in some 
situations there is a compact Lie 
group $G$ acting on $M$ and preserving $f$. 
This problem naturally arises in the $n$-body problem, where the central 
configurations are critical points of a Morse function which is invariant with 
respect to the action of $SO(n)$ (c.f. \cite{Pacella}). Another example is 
the problem of finding the number of closed geodesics of a Riemannian metric; where   
the Lie group is $S^1$ acting on the loop space of the underlying manifold and 
the Morse function is given by the energy of loops. Actually this last example was amongst 
the first applications of Morse theory, worked out 
by Morse in \cite{Morse-Calculus}; see also \cite[chapter 3]{Milnor-Morse}. In these cases,   
the connected components of critical level sets of $f$ are clearly orbits of the action. We can simply 
ignore the invariance of $f$ with respect to the action and  get lower-bounds for 
the number of these critical levels using Morse-Bott theory
(see \cite{Bott-Non} and \cite[page 344]{Bott-oldandnew}). However, 
as it is clearly explained 
in \cite[pages 351-355]{Bott-oldandnew}, to get better results 
one has to consider the $G$-invariance of $f$, and this requires an 
appropriate cohomology theory that takes account of the group action. 
This is the \emph{equivariant cohomology theory} which was introduced originally by 
H. Cartan in 1940s. 



A. G. Wasserman in \cite[section 4]{Wasserman} developed 
equivariant Morse theory for a general compact Lie group $G$, showing that one could use a $G$-invariant Morse 
function to compute equivariant cohomology. 
The Morse-Smale-Witten complex for equivariant cohomology 
is not yet constructed for a general Lie group $G$. 
The particular case $G=S^1$ 
is recently treated by M. J. Berghoff \cite{Berghoff} and provides, 
as a by-product, a proof for the 
equivariant Morse inequalities, whose precise statement may be found in 
\cite[page 149]{Wasserman} or in \cite[page 351]{Bott-oldandnew}. The methods of \cite{Berghoff} 
are apparently difficult to be generalized to a general Lie group $G$. 

In this paper we also consider the case $G=S^1$ and prove the 
equivariant Morse inequalities by adopting Witten's method to deform the Borel complex. 
Applying this method to a general Lie group (except torus groups) is not straightforward. 
Our approach will be interesting for those who 
are interested in analytical methods rather than in topological ones. 
It is also valuable for those who are interested in Morse inequalities 
rather than in the equivariant cohomology groups. Moreover, 
the analytic approach is more flexible and can be adapted to different 
situations. It has also been used to deal 
with Morse inequalities on manifolds with boundary in \cite{Zadeh-Morse}, and to establish the  
so-called delocalized Morse inequalities in \cite{Zadeh-Delocalized}. 
In \cite{WenLu} similar methods are 
used to establish equivariant Morse-Bott inequalities, where a finite group $G$ acts on a closed  
manifold $M$.

The structure of the paper is as follows. 
In section \ref{section2} we give the definition of equivariant 
cohomology theory by introducing the Borel complex. Then we give the precise 
statement of the equivariant Morse inequalities in theorem \ref{thm1}. 
In section \ref{section3} we establish the Hodge theory for the equivariant 
complex and introduce the Witten deformation of the Borel complex. 
Then we prove an infinite number of quite general equivariant analytical Morse inequalities 
in theorem  \ref{thm2}. The asymptotic behaviour of these inequalities leads finally to 
the Morse inequalities. Lemma \ref{localEu} will be used in last section  
to show that the number of Morse inequalities is actually finite. 
Lemma \ref{keylemma} reduces our problem to computing the kernel of some 
elliptic operators on Euclidian spaces. The computation of these kernels is the subject of 
section \ref{section4} and leads to the proof of the main theorem in the last 
paragraph of the paper.

\section{Equivariant cohomology and the statement of the main theorem}\label{section2}

Let $G$ be a compact Lie group that acts on a topological space $M$, and let $EG\to BG$ be a 
model for the universal principal $G$-bundle. Let $M_G$ denote the Borel construction $(M\times EG)/G$. 
The equivariant cohomology of $M$ with complex coefficients, which we denote by $H^*_{G}(M)$, is the singular cohomology of 
$M_G$ with coefficient in $\C$ 

\begin{equation}\label{eqco1}
H_{G}^*(M)=H^*(M_G,\C)
\end{equation}
If $G$ acts trivially 
then $M_G=M\times BG$ and $H_{G}^*(M)=H^*(M\times BG)=H^*(M)\otimes H^*(BG)$. 
In particular the equivariant cohomology of 
a single point is the group cohomology $H^*(BG)$. 
If the action of $G$ on $M$ is free then $EG\to M_G\to M/G$ is a fibration with 
contractible fibres which implies $H_{G}^*(M)=H^*(M/G)$. 
Actually the equivariant cohomology gives a way to combine the cohomology of the space $M$ and that 
of $G$ while taking into the account the action of $G$ on $M$, c.f. \cite{Atiyah&Bott}. 

When $M$ is a differentiable manifold and the action is 
smooth there is a more geometric approach to the construction of the equivariant cohomology.
If $G$ is the circle group $S^1$,  this construction is suitable for applying analytic methods 
and constructing an appropriate Hodge theory. 
Therefore, from now on,  we will consider this case, and denote the circle group by  $\mathbb{T}$.  
In this paper we will deal with Borel's  definition of equivariant de Rham complex as is explained,  
amongst other applications,  in \cite{Atiyah&Bott} and \cite{Bott-indomitable}.  

Let $M$ be a closed $n$-dimensional manifold which is acted upon by the group $\mathbb{T}$. 
This action is assumed to be smooth and not 
necessarily free. The Lie algebra of $\T$ is $\R$ with a fixed element $1$. 
This element generates a vector field $v$ on $M$ which is tangent to the orbits 
and vanishes at fixed points of the action. The vector field $v$ is called the infinitesimal 
generator of the action and we denote its  $t$-time flow by $\phi_t$. By convention, the 
period of $v$ and then the period of the flow $\phi_t$, is $2\pi$.  
Let $\Omega_\mathbb{T}^*(M)\subset \Omega^*(M)$ consists of all invariant differential forms 
$\omega$, i.e., those satisfying 
$\phi_t^*(\omega)=\omega$ for $t\in\R$, or equivalently $\mathcal L_v(\omega)=0$, 
where $\mathcal L_v$ is the Lie derivative with respect to $v$. 
Consider the algebra $\Omega_{eq}^*(M):=\C[t]\otimes \Omega_\mathbb{T}^*(M)$. 
If we write $t\in\mathrm{g}^*$ for dual element corresponding to $1$, i.e. $t(1)=1$, 
then this is the algebra of all polynomial functions on the Lie algebra $\mathrm{g}=\R$ 
with values in $\Omega_\mathbb{T}^*(M)$. In other words, this is the algebra of polynomials of $t$ with coefficients in $\Omega_\mathbb{T}^*(M)$. 
This algebra is graded by the rule 
 $\deg(t^i\otimes\omega)=2i+\deg(\omega)$. The linear map   

\begin{gather}
d_{eq}:\Omega_{eq}^*(M)\to \Omega_{eq}^{*+1}(M)\notag\\
d_{eq}(t^i\otimes\omega)=t^i\otimes d\omega+t^{i+1}\otimes i_v\omega\label{deq}
\end{gather}
is a differential, i.e. $d_{eq}^2=0$, and increases the degree by one. 
The \emph{equivariant de Rham cohomology groups} $H^*_{eq}(M)$ are the cohomology groups of 
this graded differential complex. 
It turns out that these groups are isomorphic to the groups $H_G^*(M)$ introduced by 
\eqref{eqco1} when $M$ is a smooth manifold \cite[page 28]{GuilStern}.
The cohomology group $H^k_{eq}(M)$ is a finite dimensional complex vector space.  
So one can define the equivariant Betti numbers by 
$\beta_{eq}^k:=\dim H^k_{eq}(M)$. If the action 
$S^1\times M\to M$ is free, then $\tilde M=M/S^1$ is a smooth manifold, and 
the equivariant cohomology groups of $M$ are canonically isomorphic to the 
de-Rham cohomology of $\tilde M$. The equivariant cohomology of a point is just 
the algebra $\C[t]$ (with $\deg t=2$), which is 
isomorphic to the de-Rham cohomology of $BS^1=\C P^\infty$.

Let $f:M\to \R$ be an invariant smooth function, i.e. $v.f=0$. An orbit $o$ is critical if 
one point on it (hence all of its points) is a critical point for $f$. For $x\in o$ 
let $N_x$ stands for the quotient space $T_xM/T_xo$. A Riemannian metric $g$ on $M$ allows us to identify 
$N_x$ with the orthogonal complement of $T_xo$. Let $\nabla$ be the Riemannian connection on $TM$ 
associated to the metric.  
For $x\in M$ and $X,Y\in T_x M$ the Hessian of $f$ is a symmetric bi-linear form 
defined as follows 

\[H_f(X,Y)=X.(Y.f)-(\nabla_XY).f\]
Because $\T$ is compact it is always possible, 
through an averaging procedure, to assume $g$ be $S^1$-invariant.
With this assumption, if $X$ or $Y$ belong to 
$T_xo$ then $H_f(X,Y)=0$. 
Therefore the Hessian defines a well defined symmetric bi-linear form on $N_x$. 
Using the Riemannian metric, we can identify $N_x$ with the orthogonal 
compliment of $T_xo$.   
We denote the restriction of $H_f$ to $N_x\subset T_xM$ by $\bar H_f$. We say a critical orbit
$o$ is transversally non-degenerate (or simply non-degenerate) if $\bar H_f$ is 
non-degenerate at $x$ , and therefore at any other point of $o$. Consequently, we may write an orthogonal direct sum $N_x=N_x^-\oplus N_x^+$ 
such that the Hessian is negative definite on $N_x^-$ and positive definite on $N_x^+$. It turnes out that these spaces can be chosen 
in a manner that form vector bundles $N$, $N^-$ and $N^+$ over $o$ such that 
\begin{equation}\label{Asd}
N=N^-\oplus N^+
\end{equation}

The Morse index of such an orbit is 
the dimension of the 
maximal subspace of $N_x$ on which the Hessian is negative-definite, i.e. the dimension of $N_x^-$.  
In what follows, we reserve the notation 
$o$ for a non-trivial orbit and we denote a trivial orbit by its geometric image, 
that is a point $p$ in $M$. If all critical points and orbits of the smooth invariant 
function $f$ are transversally non-degenerate, then $f$ is called a Morse function. 
Let $c_k$ denote the number of critical 
points of the Morse function $f$ with Morse index $k$. Let also $d_k$ denote the number of those critical orbits $o$ 
with Morse index $k$ such that the bundle $N^-$ defined in\eqref{Asd} is trivial. 
Our aim is to provide an analytic proof for the following 
equivariant Morse inequalities via Witten deformation:   

\begin{thm}\label{thm1}
With $\tilde c_k:=d_k+c_k+c_{k-2}+c_{k-4}+...$ 
the following inequalities hold for $k=0,1,2,\dots$  

\[\tilde c_k-\tilde c_{k-1}+\dots\pm \tilde c_0\geq \beta_{eq}^k-\beta^{k-1}_{eq}+
\dots\pm \beta_{eq}^{0}.\]
Actually the inequalities for $k\geq n+1$ are equivalent to the  inequality for $k=n$.
\end{thm}

Note that, when the action is free, these inequalities reduce to the ordinary 
Morse inequalities for the function 
$\tilde f\colon \tilde M\to\R$, while they reduce to the Morse inequalities 
for the function $f$ when the action is trivial. 

Equivariant Morse inequalities are mostly stated in the literature in slightly more complicated topological expression. 
Here we explain briefly how these expressions are equivalent to above inequalities  
for circle group $\T$ and for equivariant cohomolgy with coefficients in $\C$: 
The equivariant Poincar\'e series for a closed $G$-manifold $M$ is defined by 
\[P_t^G(M)=\sum_kt^j\dim H^k_{G}(M,\C)=\sum_kt^j\beta_{eq}^k.\]
Here the equivariant cohomology groups have coefficients in $\C$.  Similarly the equivariant Morse series of a Morse function
is defined as follows
\[M_t^G(f)=\sum_N t^{\lambda_N}P_t^G(N;\theta^-).\]
Here $N$ runs over the critical manifolds of $f$, and $\lambda_N$ is the Morse index of the critical manifold $N$. 
Real line bundle $\theta^-$ is the orientation bundle associated to $N^-$, and in $P_t^G(N;\theta^-)$ cohomology groups are considered with 
coefficients in flat line bundle $\theta^-$. 
The equivariant Morse inequalities as stated, e.g. in \cite[page 351]{Bott-oldandnew} assert that there is a series 
$Q_t(f)=q_0+q_1t+q_2t^2+\dots$, with positive coefficients, such that 
\begin{equation}\label{BMG}
M_t^G(f)-P_t^G(M)=(1+t)Q_t(f).
\end{equation}
Let see what information this formulation of Morse inequalities provide when $G$ is the circle group $\T$. 
In this case the equivariant Hodge isomorphism (see \eqref{had}  and remark \ref{taghi}) provide the following relation 
(see also \cite[page 351-353]{Bott-oldandnew}): 
$P_t^\T(p)=1+t+t^2+t^4+\dots$, while $P_t^\T(o,\theta^-)=0$ if $\theta^-$ is not orientable, and $P_t^\T(o,\theta^-)=1$ otherwise.
Using these equalities, relation \eqref{BMG} takes the following form 
\[\sum_{o\in O^+}t^{\lambda_o}+\sum_pt^{\lambda_p}(1+t^2+t^4+\dots)-\sum_k\beta_{eq}^k\,t^j=(1+t)(q_0+q_1t+q_2t^2+\dots)\]
where $O^+$ consists of all orbits whose line bundle $\theta^-$ is orientable. 
By comparing coefficients of powers of $t$ from both sides of this relation, and using $q_k\geq0$, 
it is easy to see that this relation implies inequalities of our theorem \ref{thm1} and vice versa.

\begin{rem}\label{dinv}
One might be interested in forming an "invariant cohomology" $H^*_{inv}(M)$ by considering 
the following complex instead of \eqref{deq}  
\[d: \Omega_\mathbb{T}^*(M)\to \Omega_\mathbb{T}^{*+1}(M)\]
Actually this complex produces the ordinary de Rham cohomology $H^*(M)$. 
To see this consider a class $[\omega]$ is $H^*(M)$. Since $\mathbb{T}$ is connected, 
each element $g\in \mathbb{T}$, as a map on $M$, is homotopic to identity map and therefore 
$[g^*\omega]=g^*[\omega]=[\omega]$ in $H^*(M)$. By averaging and normalizing the relation 
$[g^*\omega]=[\omega]$ over $\mathbb{T}$, with respect to an invariant measure, we conclude that each class in 
$H^*(M)$ is represented by a $\mathbb{T}$-invariant differential forms. Therefore the morphism 
$\phi: H^*_{inv}(M)\to H^*(M)$ defined by $\phi([\omega])=[\omega]$ is surgective. 
It is indeed injective and therefore an isomorphism. 
To see this let $\phi([\omega])=0\in H^*(M)$, then we have $\omega=d\eta$ and $\omega$ is $\mathbb{T}$-invariant. 
Then $\omega=d\, h^*\eta$ for $h\in \mathbb{T}$. 
This implies $\omega=d\alpha$, where $\alpha\in \Omega^*_{\mathbb{T}}(M)$ is the average of $h^*\eta$ on $\mathbb{T}$. 
Consequently $[\omega]=0\in H^*_{inv}(M)$ which implies the injectivity of $\phi$ and then the 
isomorphism $H_{inv}^*(M)\simeq H^*(M)$. As a result, the  
"$\mathbb{T}$-invariant Euler characteristic" and ordinary Euler characteristics are equal. 
We need this result in the following section. 
\end{rem}

 \section{Equivariant Hodge theory and analytic Morse inequalities}\label{section3}
 
 For our purposes in this paper, we need to establish an equivariant version of Hodge theory. 
 The space of $\T$-invariant 
 differential forms $\Omega^*_\mathbb{T}(M)$ is endowed with an inner product coming from the 
 Riemannian metric $g$ and its
 natural lifting to the exterior algebras $\wedge^*(T_p^*M)$ for all $p\in M$. 
 The formal dual of the exterior 
 differential $d$ is the operator $d^*$. We define a scalar product 
 on $\C[t]\otimes \Omega_\mathbb{T}^*(M)$ by the bi-linear extension of the following formula 
 
 \begin{equation}\label{eqin}
 \langle t^i\otimes\omega~,~t^j\otimes\eta\rangle=\delta_{ij}\langle\omega~,~\eta\rangle 
 \end{equation}
 It is easy to verify that the formal adjoint of the equivariant exterior derivative $d_{eq}$ 
 (see \eqref{deq})  
 with respect to the scalar product \eqref{eqin} is the linear extension of the following 
 operator of grade $-1$ 
 
 \begin{gather}
d_{eq}^*:\Omega_{eq}^*(M)\to \Omega_{eq}^{*-1}(M)\notag\\
d_{eq}^*(t^i\otimes\omega)=t^i\otimes d^*(\omega)+\delta_{i}t^{i-1}
\otimes v^*\wedge\omega\label{dseq}
\end{gather}
Here $v^*\in \Omega^1_\mathbb{T}(M)$ 
is the dual of the vector field $v$, while for $i=0,1,2,\dots$ 
\[ \delta_i=
\left\{
\begin{array}{cc}
0 &\text{if } i=0\\
1 & \text{if } i>0
\end{array}
\right.
\]
 
 The equivariant Laplacian operator is defined by 
 $\Delta_{eq}:=d_{eq}^*d_{eq}+d_{eq}d_{eq}^*$.  
 The following relation is a very direct consequence of definitions and gives 
 the action of $\Delta_{eq}$ on a term 
 like $t^i\otimes\omega$. If $\Delta=dd^*+d^*\!d$, then we have
 
 \begin{align}\label{exlap}
\Delta_{eq}(t^i\otimes\omega)=&t^i\otimes (\Delta\omega+v^*\wedge i_{v}\,
\omega+\delta_{i}\,i_{v}\,v^*\wedge\omega)\\
   &+t^{i+1}\otimes(i_v\,d^*\omega+d^*\,i_v\omega)\notag\\
   &+\delta_{i}\,t^{i-1}\otimes dv^*\wedge\omega\notag
\end{align}

 We denote the restriction of $\Delta_{eq}$ to $\Omega_{eq}^k(M)$ by $\Delta^k_{eq}$. This is a 
 second order elliptic differential operator, which is formally self-adjoint 
 with respect to the inner product \eqref{eqin}. 
 We can  
 apply the ordinary completion procedure to construct the Hilbert spaces 
 $L^2(M,\wedge_{eq}^kTM)$ or Sobolev spaces 
 $W^\alpha(M,\wedge_{eq}^kTM)$. 
  Therefore, exactly as in the classical Hodge 
 theory for the Laplacian (e.g. through the construction of heat operator as in 
 \cite[chapter 3]{Rosenberg-Laplacian}) we obtain the  
 following equivariant Hodge isomorphisms for $k=0,1,2,\dots$
  
 \begin{equation}\label{had}
H_{eq}^k(M)\simeq \ker \Delta_{eq}^k
\end{equation}

\begin{rem}\label{taghi}
Let $G=\T$ be the circle group acting trivially on a point $p$. 
This equivariant Hodge isomorphism, alongside the explicit expansion \eqref{exlap} readily imply that $H^{2k+1}_{eq}(p)=0$, while 
 $H^{2k}_{eq}(p)$ is one dimensional generated by $t^k$, for $k=0,1,2,\dots$. 
 If $\T$ acts on itself with nontrivial orbit $o$ (therefore $v$ is no-zero), then 
 equivariant Hodge isomorphism and \eqref{exlap} imply the following: If $k>0$ then $H^{k}_{eq}(o)=0$, while $H^{0}_{eq}(o)=\C$ is generated 
 by does functions of $\T$ that are invariant under the action of $\T$ (compare with \cite[page 245]{BottTu}). We could easily define a
 version of equivariant cohomology with local coefficients and prove the corresponding Hodge isomorphism to conclude that if 
 $\ell$ is a non-orientable real line bundle over $o$ then $H^k_{eq}(o,\ell)=0$ for $k=0,1,2,\dots$. 
\end{rem}

 It is clear from expansion \eqref{exlap} that for $k\geq n-1$, multiplication by $t^i$ 
 gives rise to the isomorphism  $\Omega_{eq}^k(M)\simeq\Omega_{eq}^{k+2i}(M)$ and the action 
 of the Laplacian  is linear with respect to this multiplication. Therefore,  for $k\geq n-1$ we have 
 $t^i\otimes\ker \Delta^k_{eq}=\ker \Delta_{eq}^{k+2i}$ which implies, by  \eqref{had}, 
the following results 

\begin{equation}\label{2eqper}
H_{eq}^k(M)\simeq H_{eq}^{k+2i}(M)
\end{equation} 
We use the following lemma to prove that the equivariant Morse inequalities  stop 
beyond degree $n$. 
\begin{lem}\label{localEu} 
For $k\geq n$ the following equalities hold

\[\beta_{eq}^{k}-\beta_{eq}^{k+1}=(-1)^k\chi(M)\]
where $\chi(M)$ is the Euler characteristic of $M$.
\end{lem}
\begin{pf}
Using \eqref{2eqper} it is enough to prove the lemma for $k=n$. 
It is clear from \eqref{deq} that $d_{eq}$ is $\C[t]$-linear. However,   
due to the term $\delta_{i}$ in \eqref{dseq}; the operator $d^*_{eq}$ 
is $\C[t]$-linear on $\Omega^{k}_{eq}$ only for $k\geq n$. Because 
$t\otimes\Omega_{eq}^{n-1}(M)=\Omega_{eq}^{n+1}$ and 
$t\otimes\Omega_{eq}^{n}(M)=\Omega_{eq}^{n+2}$ we can define the following grad-reversing  operator 

\begin{gather*}
\bar D_{eq}\colon \Omega_{eq}^n\oplus \Omega_{eq}^{n+1}\to \Omega_{eq}^n\oplus \Omega_{eq}^{n+1}, 
\end{gather*}
where $\Omega_{eq}^n\oplus \Omega_{eq}^{n+1}$ is given a $\Z_2$-grading, and on the even-summand, 
i.e. $\Omega_{eq}^n$ we define $\bar D_{eq}$ to be $d_{eq}+td^*_{eq}$, while on the odd-summand $\Omega_{eq}^{n+1}$ 
we define it to be 
$t^{-1}d_{eq}+d^*_{eq}$. 

The operators $\bar d_{eq}$ and $\bar d_{eq}^*$ are formal adjoint of 
each other with respect to the inner product \eqref{eqin}. 
Therefore, 
$\bar D_{eq}$ 
is grading reversing and formally self-adjoint.  
It is clear from the definitions that $\bar D_{eq}^2=\Delta_{eq}$, thus, $\bar D_{eq}$ 
is an elliptic differential operator and its Fredholm index is given by 
 
 \[\ind \bar D_{eq}=\dim\ker \Delta_{eq}^n-\dim\ker \Delta_{eq}^{n+1}=
 \beta_{eq}^n-\beta_{eq}^{n+1}.\]
On the other hand, as differential operators 
on $\Omega_{eq}^n\oplus \Omega_{eq}^{n+1}$, the elliptic operators  
$\bar D_{eq}$ and the classical de Rham operator $D:=d+d^*$, acting on the 
$\mathbb{T}$-invariant differential forms, have the same principal symbols. 
Therefore they are homotop in the space of Fredholm operators.  
Because Fredholm index 
is homotopy invariant, we get the equality \(\ind \bar D_{eq}= \ind D\). The right side of this 
equality could be named as $\mathbb{T}$-invariant Euler characteristic, which is actually 
equal to the ordinary Euler characteristic 
of $M$ by remark \ref{dinv},  and this completes the proof of the lemma. 
\end{pf}

We need a deformed version of the equivariant Laplacian that we introduce here.  
Given a positive parameter $s$ and using the Morse function $f$,  we define the the Witten 
deformation of $d_{eq}$  as the operator 
$d_{eq,s}:=e^{-sf}d_{eq}e^{sf}=d_{eq}+sdf\wedge$. Its formal dual is 
$d_{eq,s}^*=d_{eq}^*+s\nabla_f\lrcorner$, where $\nabla_f$ stands for the gradient of $f$. 
The associated deformed equivariant Lapacian 
$\Delta_{eq,s}:=d_{eq,s}^*d_{eq,s}+d_{eq,s}d_{eq,s}^*$ has the following expansion, c.f. 
\cite[Lemma 9.17]{Roe-elliptic}

\begin{equation}\label{deflap}
\Delta_{eq,s}=\Delta_{eq}+s^2|df|^2+s\mathrm{H}_f
\end{equation}
Here $\mathrm{H}_f$ is the following operator, where $\{e_i\}_i$ is a local orthonormal 
base for $TM$ and $L_{e_i}$ and 
$R_{e_i}$ are respectively the left and right Clifford multiplication by $e_i$ 
(see \cite[page 126]{Roe-elliptic}) 

\begin{equation}\label{clhes}
\mathrm{H}_f=\sum_{i,j} H_f(e_i,e_j)L_{e_i}R_{e_j}
\end{equation}

It is clear that $\Omega_{eq}^*(M)$ with the differential $d_{eq,s}$ is a 
graded differential complex, 
so we may define the deformed cohomology spaces $H_{eq,s}^k(M)$. Nevertheless, these 
cohomology spaces are not new objects,   
and multiplication by $e^{-sf}$ provides the following isomorphism 

\begin{equation}\label{siso}
 H_{eq}^k(M)\simeq H_{eq,s}^k(M). 
\end{equation}

Using the expansion \eqref{deflap} and Duhamel's formula \cite[page 97]{Rosenberg-Laplacian}, 
one is able to construct the heat 
operator associated to $\Delta_{eq,s}$ which establishes a Hodge theory 
and provides the isomorphisms 
$\ker \Delta_{eq,s}^k\simeq H_{eq,s}^k(M)$. This and \eqref{siso} give the following relation

\begin{equation}\label{betsbet}
\beta_{eq}^k=\dim(\ker \Delta_{eq,s}^k)
\end{equation}
 Using \eqref{deflap}, it is clear that $\Delta_{eq,s}^k$ is a non-negative and second order elliptic 
 differential operator on 
 $L^2(M,\wedge_{eq}^kTM)$. Therefore,  
 given a smooth and rapidly-decreasing positive function $\phi$ on $\R^{\geq0}$ satisfying $\phi(0)=1$, 
 the operator $\phi(\Delta_{eq,s}^k)$, being a smoothing operator on $L^2(M,\wedge_{eq}^kTM)$ , 
 is a trace class operator, and we denote its trace by 
 
 \begin{equation*}
  \mu^k_{eq,s}=\tr\phi(\Delta_{eq,s}^k)~~;\hspace{1cm} k=0,1,...,n
 \end{equation*}
  The following equivariant analytic Morse inequalities are our departure to our proof for    
 theorem \ref{thm1} (see \cite{Roe-elliptic} for the non-equivariant version) 
 
 \begin{thm}[the analytic equivariant Morse inequalities]\label{thm2}
With the above notations, the following inequalities hold for $k=0,1,2,\dots$
\[\mu^k_{eq,s}-\mu^{k-1}_{eq,s}+\dots\pm \mu^0_{eq,s}\geq
 \beta_{eq}^k-\beta_{eq}^{k-1}+\dots\pm\beta_{eq}^0\]
 \end{thm}
\begin{pf}
If we put $\beta_{eq,s}^k=\dim(\ker \Delta_{eq,s}^k)$ then by \eqref{betsbet} 
the above inequalities 
are equivalent to the followings 
\[\mu^k_{eq,s}-\mu^{k-1}_{eq,s}+\dots\pm \mu^0_{eq,s}\geq
\beta_{eq,s}^k-\beta_{eq,s}^{k-1}+\dots\pm\beta_{eq,s}^0\]
The argument in the proof of the proposition 14.3 of \cite{Roe-elliptic} applies 
to the deformed Laplacian and gives these inequalities. For the sake of 
completeness we give a very brief account of this proof.  
The spectrum of the deformed Laplacian $\Delta_{eq,s}^k$ is discrete, so there is a 
rapidly decreasing function $\tilde\phi$ on $\R$ which vanishes on non-zero elements 
of the spectrum such that $\tilde{\phi}(0)=1$. Therefore  
$\beta_{eq,s}^k=\tr\tilde{\phi}(\Delta_{eq,s}^k)$ which implies 
$\mu^k_{eq,s}-\beta_{eq,s}^k=\tr((\phi-\tilde{\phi})\Delta_{eq,s}^k)$. The relation 
$(\phi-\tilde{\phi})(x)=x\psi(x)^2$ defines a rapidly decreasing function $\psi$ on $\R$,  
and one get the following 
\begin{equation}\label{zozi}
\mu^k_{eq,s}-\beta_{eq,s}^k=\tr \Delta_{eq,s}^k\psi(\Delta_{eq,s}^k)^2
\end{equation}
Let $H_j$ denote the $L^2$-Hilbert space generated by $\Omega_{eq}^j(M)$. 
Using the relation $\Delta_{eq,s}^j=d_{eq,s}^*d_{eq,s}+d_{eq,s}d_{eq,s}^*$ 
one gets easily the following relation 

\[\tr\{d_{eq,s}d_{eq,s}^*\psi(\Delta_{eq,s}^j)^2\}_{|H_j}=
\tr\{d_{eq,s}^*d_{eq,s}\psi(\Delta_{eq,s}^{j-1})^2\}_{|H_{j-1}}\]
Using this with an alternating summation from $j=k$ to $j=0$ on equation \eqref{zozi} we get 
the following relation 

\[(\mu^k_{eq,s}-\beta_{eq,s}^k) -(\mu^{k-1}_{eq,s}-\beta_{eq,s}^{k-1})+\dots\pm (\mu^0_{eq,s}-\beta_{eq,s}^) =
\tr\{d_{eq,s}^*d_{eq,s}\psi(\Delta_{eq,s}^k)^2\}_{|H_{k}}\]
Since $d_{eq,s}^*d_{eq,s}\psi(\Delta_{eq,s}^k)^2$ is a non-negative operator, 
the right side of the above relation is non-negative in general and this 
gives the equivariant Morse inequalities. 
\end{pf}

To prove theorem \ref{thm1}, we will study the asymptotic behavior of $\mu^k_{eq,s}$ when $s$ 
goes toward infinity and apply theorem \ref{thm2}. 
Since $\phi$ is rapidly deceasing, the operator 
$\phi(\Delta^k_{eq,s})$ is smoothing and 
has a smooth kernel 
\[\phi(\Delta^k_{eq,s})\omega\,(p)=\int_MK_{s}^k(p,q)\,\omega(q)\,d\mu_g(q)\]
Here $K_{s}^k(p,q)$ is an element of $\wedge_{eq}^kT_pM\otimes\wedge_{eq}^kT_q^*M$ and 
$\mu_g$ is the Riemannian volume element associated to the invariant metric $g$. 
Therefore for $k=0,1,2,\dots$ we have 

\begin{equation}\label{strac}
 \mu^k_{eq,s}=\int_M\tr K_{s}^k(p,p)\,d\mu_g(p)
\end{equation}

\noindent We recall the following relation from \eqref{deflap} 
\[\Delta_{eq,s}=\Delta_{eq}+s^2|df|^2+s\mathrm{H}_f\]
Let's restrict ourself to a complement set of an open neighbourhood of the critical levels  
(i.e. the union of the critical points and the critical orbits) 
of $f$, where $|df|\geq c>0$.  
Here $\Delta_{eq}$ is non-negative, while the term $s^2|df|^2$  dominates the term $s\mathrm H_f$ when $s$ 
goes to infinity . Therefore, informally speaking,  on the smooth sections supported in this set, operator 
$\Delta_{eq,s}^k$ get bigger and bigger t when $s$ goes to infinity. Consequently on this set $\phi(\Delta_{eq,s}^k)$ goes to zero  
when $s$ goes to infinity. 
This argument can actually provide a rigorous proof for the 
following lemma by using 
finite propagation speed property of the wave operator and the Friedrich extension theorem. This is done in the proof 
of the non-equivariant case in \cite[Lemma 14.6]{Roe-elliptic} and can be literally applied to our equivariant 
context to provide a proof for the following lemma: 

\begin{lem}\label{tracout}
When $s$ goes toward infinity, the smoothing kernel $K_{s}^k(p,q)$ goes uniformly to zero 
when $p$ or $q$ belong to a complement of an open neighbourhood of the critical 
levels of $f$. 
\end{lem}

For $\rho>0$ let $N_{4\rho}(p)$ and $N_{4\rho}(o)$ denote, respectively, 
the $4\rho$-neighbourhood of the critical point $p$ and the critical orbit $o$. 
Let also $\phi_p$ and $\phi_{o}$ denote equivariant non-negative  
smooth functions on $M$  which are supported, respectively, 
in $N_{3\rho}(p)$ and in $N_{3\rho}(o)$ such that $\phi_p=1$ on $N_{\rho}(p)$ and 
$\phi_{o}=1$ on $N_\rho(o)$. Point-wise multiplication by these functions 
defines operators on equivariant 
differential forms. The following corollary comes up as a 
very direct result of the above lemma:

\begin{cor}\label{cor1}
The following relation holds 

\[\lim_{s\to\infty} \mu^k_{eq,s}=\lim_{s\to\infty}\tr\phi(\Delta_{eq,s}^k)=
\sum_p \lim_{s\to\infty}\tr(\phi_p\phi(\Delta_{eq,s}^k))+
\sum_{o}\lim_{s\to\infty} \tr(\phi_{o}\phi(\Delta_{eq,s}^k)\]
\end{cor}

Let $\mathbb B_{a}^n(0)$ denote the ball in $\R^n$ with center $0$ and radius $a$. 
By choosing $\rho$ sufficiently small and using an equivariant version of partition of unity, 
we can assume that $N_{4\rho}(p)$  and $\mathbb B_{4\rho}^n(0)$ are isometric 
and in this isometry the point $p$ corresponds to $0$. This is also true for 
$N_{4\rho}(o)$ and $S^1\times \mathbb B_{4\rho}^{n-1}(0)$ where $o$ corresponds 
to $S^1\times \{0\}$, provided that $N_{4\rho}(o)$ be orientable. 
We will consider the non orientable case later. 
Let $L_s^k$ and $\bar{L}_s^k$ denote, respectively, differential operators 
on  $\Omega_{eq}^k(\R^n)$ and $\Omega_{eq}^k(S^1\times\R^{n-1})$ such that 
with respect to above isometries 

\[\Delta^k_{eq,s}|_{N_{4\rho(p)}}=L^k_{s}|_{\mathbb B_{4\rho}^n(0)}\hspace{0.5cm} \text{ and } 
\hspace{0.5cm}\Delta^k_{eq,s}|_{N_{4\rho}(o)}=\bar{L}^k_{s}|_{\mathbb S^1\times B_{4\rho}^{n-1}(0)}\] 

\noindent Then, through a standard argument, based on Fourier inversion formula and 
finite propagation speed of wave operators, the following equalities hold 

\[\phi(\Delta_{eq,s}^k)(\omega_1)=\phi(L_{s}^k)(\omega_1)\hspace{0.5cm} \text{ and } 
\hspace{0.5cm} \phi(\Delta_{eq,s}^k)(\omega_2)=\phi(\bar{L}_{s}^k)(\omega_2)\]
provided that the Fourier transform $\hat \phi$ of $\phi$ is supported in 
$(-\rho~,~\rho)$ and the support of $\omega_1$ and $\omega_2$ are included, respectively, 
in $\mathbb B_{3\rho}^n(0)$ and $ S^1\times\mathbb B_{3\rho}^{n-1}(0)$. Therefore,

\[\tr(\phi_p\phi(\Delta_{eq,s}^k))=\tr(\phi_p\phi(L_s^k))\hspace{0.5cm} \text{ and } 
\hspace{0.5cm} \tr(\phi_o\phi(\Delta_{eq,s}^k))=\tr(\phi_o\phi(\bar{L}_s^k))\]
These equalities and corollary \ref{cor1} together  prove the following lemma 

\begin{lem}\label{keylemma}
Provided that the support of $\hat\phi$, the Fourier transform of $\phi$,  
is included in $(-\rho,\rho)$, and with above notations,
the following relation holds 
\begin{equation}\label{sasa}
\lim_{s\to\infty} \mu_{eq,s}^k=\lim_{s\to\infty} \tr\phi(\Delta_{eq,s}^k)=
\sum_p \lim_{s\to\infty}\tr(\phi_p\phi(L_s^k))+
\sum_{o}\lim_{s\to\infty} \tr(\phi_{o}\phi(\bar{L}_s^k))
\end{equation} 
Here $p$ runs over all critical points of $f$ while $o$ 
runs over all critical orbits of $f$ and $L_s^k$ and $\bar{L}_s^k$ are 
the local representation of $\Delta_s^k$ around $4\rho$-neighbourhood of, 
respectively, critical points and critical orbits. 
\end{lem}
In the following section we will compute the values of each term in the right side of the above relation.

\section{Localization on critical levels}\label{section4}

At the beginning of this section, let start with some well-known facts which will be used afterward.  

\begin{lem}\label{zeplo} 
Let $Z_i := [dx_i\wedge~,~dx_i\lrcorner] $ be considered as a linear map on exterior algebra generated by $dx_i$'s.  
This linear map is diagonalizable and  an element of the form 
$dx_{i_1}\wedge dx_{i_2}\wedge\dots \wedge dx_{i_j}$ is an eigenvector corresponding to eigenvalue 
$1$ if $i=i_\ell$ for one $\ell$, and to eigenvalues $-1$ otherwise. 
\end{lem}

In order to evaluate the right side of \eqref{sasa}, and therefore compute the asymptotic 
behavior of $\mu_s^k$ at the 
left side of that equality, we   need to review some spectral properties of 
\emph{harmonic oscillator operator}  (see \cite[]{Roe-elliptic})

\begin{equation}\label{harmos}
-(\frac{\partial^2}{\partial x_1^2}+\frac{\partial^2}{\partial x_2^2}+\dots+\frac{\partial^2}{\partial x_n^2})+
a^2(x_1^2+x_2^2+\dots+x_n^2)~~;\hspace{1cm} a>0
\end{equation} 
This is an unbounded self-adjoint operator acting on $L^2(\R^n)$,  the completion of $C^\infty_0(\R^n)$ with respect to 
$L^2$-norm, and provides a spectral resolution for this Hilbert space.  
The eigenvalues of this operator are $a(n+2p)$ with $p=0,1,2,\dots$. 
The eigenvector corresponding to the minimal eigenvalue  $na$ is the following function 
\begin{equation}\label{sahand}
 u_0(x):=(a\pi^{-2})^{n/4}\exp(-ax^2/2)
\end{equation}
Here $x=(x_1,x_2,\dots,x_n)$ and $x^2=x_1^2+x_2^2+\dots+x_n^2$. 
Moreover, given a compactly supported smooth function $\beta$ on $\R^n$ such that $\beta(0)=1$, then 

\begin{equation}\label{locin}
\lim_{a\to\infty}\langle\beta(x)u_0(x)~;~u_0(x)\rangle=\beta(0)=1
\end{equation}

To compute the contribution of critical levels in \eqref{sasa},, we need to have a good 
representation of the deformed equivariant 
Laplacian operators around them. This is provided by an  
equivariant version of the Morse lemma that we are going to explain. 
Let us begin with an 
equivariant version of the tubular neighbourhood theorem. Suppose that $G$ is a compact Lie group acting 
on the closed manifold $M$ and  $g$ is a Riemannian metric which is invariant under the action.
The map $\eta_x:G/G_x\to G.x$ given by $\eta_x([h])=h.x$ is a diffeomorphism, where $G.x$ and $G_x$ are the 
orbit and stablizer of $x\in M$. 
For $h\in G_x$, the derivative $T_xh\colon T_xM\to T_xM$ is an isometry that 
maps $T_xG.x$ into itself. Therefore, it induces a linear isometry 
$\phi(h):N_x\to N_x$ where $N_x\subset T_xM$ is the 
orthogonal complement of $T_xG.x$. 
In other words one has an orthonormal representation   
$\phi:G_x\to O(N_x)$. The subgroup $G_x$ has a free action on $G\times N_x$ given by 
$h.(h',v)=(h'h^{-1},\phi(h)(v))$. 
The quotient space is a vector bundle $\pi:N\to G.x$ whose fibres are isometric to $N_x$ 
(here we have used the identification $G.x= G/G_x$). 
The action of $G$ on $G\times N_x$ is given by $h.(h',v)=(hh',v)$, and it  
commutes with the action of $G_x$. Therefore it induces a bundle map on the 
vector bundle $N\to G.x$. The equivariant tubular neighbourhood theorem \cite[page 15]{Audin-Torus} 
asserts that there is an invariant neighbourhood $W$ of $G/G_x$, as the zero 
section of the bundle $N$, and an invariant neighbourhood $U$ of the orbit $G.x$, and an 
equivariant diffeomorphism  $\bar \eta\colon W\to U$ that extends the orbit map $\eta$ and 
makes the following diagram commutative 

\begin{equation*}
\begin{CD}  
 G/G_x @> \eta >>  G.x\\
 @V i VV         @V i VV\\
 W\subset N @> \bar \eta >>  U\subset M
     \end{CD}
\end{equation*} 
Note that by this tubular neighbourhood theorem, each fiber $N_y$ has an inner product and the action 
$g:N_y\to N_{g.y}$ is a linear isometry. Therefore $N\to G.x$ is a Riemannian vector bundle, and we have 
\emph{equivariant Morse Lemma}  \cite[lemma 4.1]{Wasserman}: For $r>0$ put 
\[N(r)=\{v_y\in N|\|v_y\|_y<r,~y\in G.x\}\]
and let $G.x$ be a non-degenerate critical manifold for an invariant function $f$. 
There is an equivariant diffeomorphism $\psi: N(r)\to U$, for some $r>0$, such that for $v\in N(r)$ we have 
$f(\psi(v))=\|Pv\|^2-\|(1-P)v\|^2$, where $P$ is an equivariant orthogonal bundle projection. 
 
Now we go back to our case, where the Lie group $G$ is the circle group $\T$ and compute the contribution of critical 
points and critical orbits in the right side of \eqref{sasa}. \\

 \noindent{\bf Contribution of  critical fixed points:} let $p$ be a fixed point 
of the action and a critical point for the Morse function $f$. We may assume that $p=0$ and $f(p)=0$. 
The equivariant tubular neighbourhood 
theorem and the Morse lemma provide a coordinate system around $p$ with respect 
to which the Riemannian metric takes the form 
\begin{equation}\label{metp}
g=dx_1^2+\dots+ dx_n^2~,
\end{equation}
the elements of $G$ act as 
elements of $SO(n)$, and $f$ takes the following form, where $m$ is the Morse index of $p$ 
\[f(x)=- x_1^2-\dots-x_m^2+ x_{m+1}^2+\dots + x_n^2.\] 
Because the elements of $\T$ are linear functions that preserve both the Euclidean norm 
$x_1^2+\dots+x_m^2+ x_{m+1}^2+\dots + x_n^2$ and the quadratic form 
$- x_1^2-\dots-x_m^2+ x_{m+1}^2+\dots + x_n^2$, they preserves the subspaces 
$(x_1,\dots,x_k,0,\dots,0)$ as well as $(0,\dots,0,x_{k+1},\dots,x_n)$, and on these subspaces the action 
preserves the Euclidean norm. By standard results from representation theory 
(or canonical forms of orthogonal operators) there are orthonormal basis for these subspaces with 
respect to which the action of $\T$, and the Morse function $f$ take the following forms 
(as in this new coordinates the form of 
$g$ and $f$ do not change, wee keep to denote this last coordinates by $x=(x_1,x_2,\dots, x_n)$) 
\begin{equation}\label{actp}
e^{i\theta}(x_1+ix_2,\dots,x_n)=
(e^{im_1\theta}(x_1+ix_2),\dots,
e^{im_q\theta}(x_{2q-1}+ix_{2q}),x_{2q+1},\dots, x_n)
\end{equation}
where $m_i \in \N$, and

 \begin{equation}\label{funp}
 f(x_1,x_2,\dots,x_{2q+1},\dots,x_n)=\epsilon_1(x_1^2+x_2^2)+\dots+
 \epsilon_q(x_{2q-1}^2+x_{2q}^2)+\lambda_{2q+1}x^2_{2q+1}+\dots+\lambda_n x_n^2
 \end{equation}
 where $\epsilon_j$'s and $\lambda_j$'s are equal to $\pm 1$ and the 
 total number of occurrence of $-1$ is equal to the Morse index of $p$. In \eqref{actp} the expressions like 
 $e^{im_1\theta}(x_{1}+ix_{2})$ denotes the multiplication of the complex numbers $e^{im_1\theta}$ 
 and $x_{1}+ix_{2}$. The vector field $v$ and its dual with respect to $g$ take 
the following form in this coordinate system 

\begin{gather}
 v=(-m_1x_2,m_1x_1,\dots,-m_qx_{2q},m_qx_{2q-1},0,\dots,0),\label{vecp}\\
 v^*=-m_1x_2dx_1+m_1x_1dx_2-\dots-m_qx_{2q}dx_{2q-1}+m_qx_{2q-1}dx_{2q} \label{vec*p}
\end{gather}

\noindent Also the Clifford hessian of \eqref{clhes} takes the following form (see \cite[page 126]{Roe-elliptic})

\begin{equation}\label{zip}
\mathrm H_f=\sum_{i=1}^n\lambda_i Z_i~;~~Z_i=[dx_i\wedge~,~dx_i\lrcorner]
\end{equation}
where $\lambda_i = \pm 1$ is the coefficient of $x_i$ in \eqref{funp}.

Using \eqref{deflap} and \eqref{zip}, the deformed Laplacian $\Delta_{eq,s}(M)$ coincides with the following operator in a small neighborhood of $p$. 

\begin{gather}
L_s\colon\Omega_{eq}(\R^n)\to \Omega_{eq}(\R^n)\notag\\
L_s (a_{j,I} \: t^j\otimes dx^I) = ((\Delta  + 4s^2|x|^2 + sC_I)a_{j,I} )  \: t^j\otimes dx^I + a_{j,I}\:M(t^j \otimes dx^I) +a_{j,I} \: K(t^j \otimes dx^I)
\label{lalab}
\end{gather}
where $\Delta$ is euclidean Laplacian in $\R^n$, while $M$ and $K$ are tensorial operators given by
\begin{gather}
 M( t^j \otimes dx^I) = t^j \otimes (v^*\wedge i_v(dx^I) + \delta_{j}i_v(v^* \wedge dx^I))\\
K(t^j \otimes dx^I) =  t^{k+1} \otimes (i_v(d^* (dx^I)) + d^*(i_v(dx^I))) + \delta_{j}t^{k-1} \otimes dv^*\wedge dx^I
\end{gather}
Here $C_I$'s are constants defined by the last equality in the following expression for $\mathrm H_f$  (the first equality is coming 
from lemma \ref{zeplo})
\begin{equation}
\mathrm H_f(t^j \otimes dx^I) = 2(\sum_{i\in I} \lambda_i -\sum_{j \notin I}\lambda_j) t^j \otimes dx^I= C_I \: t^j \otimes dx^I  \label{hesp}
\end{equation}

Note that $M$ and $K$ are nonnegative self adjoint operators which are not dependent on $s$ and the eigenvectors of 
$T_s = L_s - M - K$ are $u^l_I  \: t^j \otimes dx^I$ where $u^l_I$ is $l$-th eigenvector of the following operator
\begin{equation}
\Delta + 4s^2|x|^2+sC_I:L^2(\R^n) \to L^2(\R^n),
\end{equation}
with corresponding eigenvalue $s(4l+2n+C_I)$, independent of $j$ (see the discussion at the beginning of this section on harmonic 
oscillator operator).   
By definition of $C_I$ in \eqref{hesp}, it is clear that $C_I \geq -2n$ and the equality holds for arbitrary $j$, but exactly one $\tilde I$ 
consisting of those indices  $i$ such that the 
coefficient of $x_i^2$ in expression \eqref{funp} for $f$ is $-1$. Therefore, all nonzero element of the spectrum of $T_s$ 
are greater than or equal to $2s$, while $0$ is an eigenvalue of $T_s$ with eigenvectors $u_0  \: t^j \otimes dx^{\tilde I}$, where 
\begin{equation}
u_0 = \left(2\pi^{-2}s\right)^{n/4}e^{-s|x|^2}\label{rezmog}
\end{equation}
is introduced by \eqref{sahand}.
Let $L^k_s$, $T^k_s$, $M^k$ and $K^k$ be the restrictions of $L_s$, $T_s$, $M$ and $K$ to equivariant space 
$\Omega^k_{eq}(\R^n)$ and denote the $l$-th eigenvalue of $L^k_s$ by $\lambda^l_{k,s}$. By above discussion $T^k_s$ 
has only one zero eigenvalue if $k = 2j+|\tilde I|=2j+m$ (m is the Mosre index of $f$ at $p$) and  the other eigenvalues are greater than $s$. For such $k$ we have
\begin{align*}
\lambda^0_{k,s} = \min \limits_{|w| = 1} (L^k_sw,w) &\leq \langle L^k_s(u_0  \: t^j \otimes dx^{\tilde I}),u_0  \: t^j \otimes dx^{\tilde I}\rangle  \\ 
&= \langle(T^k_s+M^k+K^k)(u_0  \: t^j \otimes dx^{\tilde I}), u_0  \: t^j \otimes dx^{\tilde I} \rangle.
\end{align*}
Since $T^k_s(u_0  \: t^j \otimes dx^{\tilde I})=0$ and $K^k$ changes the power of $t$, we have
\begin{align*}
\langle(T^k_s+M^k+K^k)(u_0  \: t^i \otimes dx^{\tilde I}), u_0  \: t^j \otimes dx^{\tilde I} \rangle 
&= \langle M^k(u_0  \: t^j \otimes dx^{\tilde I}),u_0  \: t^j \otimes dx^{\tilde I}\rangle  \\
\end{align*}
By the definition of $M$ and the expressions \eqref{vecp} and \eqref{vec*p} whit $C:= \max \{ m_1^2 , \dots , m_q^2 \}$ we have
\begin{align*}
\langle M^k(u_0  \: t^j \otimes dx^{\tilde I}),u_0  \: t^j \otimes dx^{\tilde I}\rangle 
&\leq C \langle |x|^2u_0  \: t^j \otimes dx^{\tilde I},u_0  \: t^j \otimes dx^{\tilde I}\rangle \\
&= {{C}\over{2s\pi^n}}\int_{\R^n}|y|^2e^{-|y|^2}\,dy
\end{align*}
Therefore, if $s \to \infty$ then $\lambda^0_{k,s} \to 0$,  and consequently  $L^k_s(u_0  \: t^j \otimes dx^{\tilde I})\to 0$.

\noindent For the $l$-th eigenvalue, where $l\geq 2$
\begin{equation}
\lambda^l_{k,s} = 
\min \limits_{\substack{V \subset \Omega^k_{eq}(\R^n) \\ dim(V)=l} }
\max \limits_{\substack{w \in V \\  |w| = 1} }(L^k_s w,w) 
\geq 
\min \limits_{\substack{V \subset \Omega^k_{eq}(\R^n) \\ dim(V)=l} }
\max \limits_{\substack{w \in V \\  |w| = 1} }((T^k_s w,w) \geq 2s
\end{equation}
\newline 
Summarizing,  
if $k-m$ is a non-negative even integer, when $s$ 
goes to infinity, the smallest  eigenvalue of $L_s^k$ (which is non negative
with multiplicity one) goes toward zero, while the other elements of its spectrum go to infinity. Therefore $\phi(L_s^k)$ converges 
to the orthogonal projection on the linear space generated by  $u_0\: t^j \otimes dx^{\tilde I}$, when $s$ goes to infinity. 
Using \eqref{locin}, this proves the following lemma 

\begin{lem} \label{consump}
Let $m$ be the Morse index of $f$ at a critical point $p$. The following relation holds

\[\lim_{s\to\infty}\tr \phi_p\phi(L_s^k)=\left\lbrace
\begin{array}{cc}
1 & \text{for }m=k, k-2, k-4,\dots\\
0 & \text{otherwise}
\end{array}
\right.\]
therefore, 
\[\lim_{s\to\infty}\sum_{p}\tr\,\phi_p\phi(L_s^k)=c_k+c_{k-2}+c_{k-4}+\dots\]
where $p$ runs over the critical points of the Morse function $f$, and $c_k$ 
denotes the number of critical points with Morse index $k$.  
\end{lem}

\noindent {\bf Contribution of critical orbits:}  Let $x$ be a point on a critical orbit $o$ and $N_x$ be the orthogonal complement of 
$T_xo$.  In this case $G_x$ is a discrete subgroup $\Z_q$ of $\T$ with an orthogonal action on $N_x$.  
The action of $\T$ on $N=\T\times N_x/\Z_q$  is induced by the translation action $e^{i\psi} (\theta, w)=(\psi+\theta, w)$ on the covering $\T\times N_x$ 
and the covering map $\pi\colon \T\times N_x\to N$ is equivariant with respect to this actions. 
By the equivariant tubular neighbourhood theorem and Morse lemma there are invariant neighbourhood $U$ of 
$o$ and $N(r)$ of the zero section of the vector bundle $N $, and there is an equivariant diffeomorphism 
$\psi : N(r) \to U$, such that $f(\psi(w))=\|Pw\|^2-\|(1-P)w\|^2$, where $P$ is an equivariant orthogonal bundle projection. The 
pull back of the bundle projection $P$ by $\pi$ is a trivial bundle projection $Q$ on $N_x$ 
(considered as a bundle map on $\T\times N_x$), and  

\[(\psi  \circ \pi)^* \!f \:(\theta,w) = |Q(w)|^2 - |(I-Q)(w)|^2\]

Therefore, there is a coordinate system $(\theta,x_1,\dots ,x_{n-1})$ for the covering space $\T\times N_x$ in which the metric has the following form 
\begin{equation}\label{abc} 
g=\kappa^2d\theta^2+dx_1^2+dx_2^2+\dots+dx_{n-1}^2,
\end{equation} 
where $\kappa^2$ is a constant that will take large values in forthcoming  discussion, and 
the function $f$ is given by the following expression  
 \begin{equation}\label{funo}
  f(\theta,x_1,\dots, x_{n-1})=\epsilon_1x_1^2 + \dots + \epsilon_{n-1}x_{n-1}^2.
 \end{equation}
Moreover, the action has the following trivial representation  

\begin{equation}\label{acto}
e^{i\psi}(\theta,x_1,\dots,x_{n-1})=(\psi+\theta, x_1 , \dots , x_{n-1}).
\end{equation} 
With respect to these coordinates, Clifford hessian of \eqref{clhes} takes the following form (see \cite[page 126]{Roe-elliptic})

\begin{equation}\label{zipo}
\mathrm H_f=\sum_{i=1}^{n-1}\epsilon_i Z_i~;~~Z_i=[dx_i\wedge~,~dx_i\lrcorner]
\end{equation}

\noindent Let denote by $L_s$ the restriction of $\Delta_{eq,s}$ to $U$ and denote by $\tilde L_s$ is lifting (as a differential operator) to 
a neighbourhood of $\T$ (as the zero section) in the trivial bundle $\T\times N_x$ with coordinates $(\theta, x_1,\dots,x_{n-1})$. 
The $\T$ invariant eigenvectors of $L_s$ and $\tilde L_s$ are independent of the variable $\theta$. 
It is clear that eigenvectors of $L_s$ can be lifted to eigenvectors of $\tilde L_s$ with the same eigenvalues. 
On the other hand, $\T$-invariant eigenvectors  of $\tilde L_s$ are actually smooth functions on $\R^{n-1}$, and 
among them, those that are invariant under the action of $\Z_q$  are lifted eigenvector of $L_s$ (with the same eigenvalue). 
Therefore we need to study the spectral properties of 

\[\tilde L_s \colon\Omega_{eq}(\T \times \R^{n-1})\to \Omega_{eq}((\T \times \R^{n-1})\]

\noindent Since all invariant differential forms in $\T \times \R^{n-1}$ are independent of 
$\theta$ and $v = (1,0,\dots,0)$ and $v^* = \kappa^2 d\theta$ we have

\begin{gather}
dv^*  \omega = 0\notag \\
(i_v \circ d^* + d^* \circ i_v) \omega = (dv^* \wedge)^* \omega = 0 \notag\\
(v^* \wedge i_v + \delta_{j} i_v(v^* \wedge))(t^j \otimes d\theta \wedge dx^I) = \kappa^2 d\theta \wedge dx^I\notag \\
(v^* \wedge i_v + \delta_{j} i_v(v^* \wedge))(t^j\otimes dx^I) = \kappa^2\delta_{j} t^j\otimes  dx^I \notag\\
\end{gather}
Therefore, by lemma \ref{zeplo} and relations \eqref{funp}, \eqref{abc} and \eqref{funo} we get following expressions: 
\begin{gather}
\mathrm H_f (d\theta \wedge dx^I)  = 2 \left ( \sum_{i \in I}  \epsilon_i - \sum_{j \notin I }\epsilon_j \right ) 
d\theta \wedge dx^I = C_I d\theta \wedge dx^I \label{cia1}\\ 
\mathrm H_f ( dx^I)  = 2 \left ( \sum_{i \in I}  \epsilon_i - \sum_{j \notin I }\epsilon_j \right ) dx^I = C_Idx^I, \notag 
\end{gather}
\begin{align*}
\tilde L_s(a_{j,I} \: t^j \otimes d\theta \wedge dx^I) &= ((\Delta + 4s^2|x|^2 + sC_I + \kappa^2) a_{j,I}) \: t^j \otimes d\theta \wedge dx^I   \\
\tilde L_s(b_{j,I} \: t^j \otimes  dx^I) &= ((\Delta + 4s^2|x|^2 + sC_I + \delta_{j}\kappa^2) b_{j,I}) \: t^j \otimes dx^I  
\end{align*}
Similar to the case of critical points, eigenvectors of $\tilde L_s$ are $u^l_I\,  \: t^j \otimes d\theta \wedge dx^I$ and $u^l_I  \: t^j\otimes  dx^I$, 
where $u^l_I$ is $l$-th eigenvector of the following operator, corresponding to the eigenvalue $s(4l+2n-2+C_I)$
\begin{equation*}
\Delta + 4s^2|x|^2+sC_I:L^2(\R^{n-1}) \to L^2(\R^{n-1})
\end{equation*}
By definition of $C_I$ in \eqref{cia1}, it is clear that $C_I \geq -2n+2$ and the equality holds for
 $\tilde I$, consisting of those indices  $i$ such that the 
coefficient of $x_i^2$ in 
 expression \eqref{funo} for $f$ is $-1$. Therefore, all  non zero eigenvalues of $\tilde L_s$ go to infinity when $s$ and $\kappa$ go to 
 infinity, while its kernel is $\tilde w_0 = u_0  \: t^0 \otimes dx^{\tilde I}$, where 
\begin{equation}
u_0 = \left(\pi^{-1}\sqrt {2s}\right)^{n/2}e^{-s|x|^2}\label{rezmog}.
\end{equation}
Since $dx^{\tilde I} $ is the volume form of $\ker Q = \pi^*\ker P$, 
it is invariant under the action of $\Z_q$ if and only if $\ker P=N^-$ is orientable (see \eqref{Asd} for definition of $N^-$). 
Of course $u_0$ is  invariant under all orthogonal transformations of $\R^{n-1}$. Therefore, $\tilde w_0$ is $\Z_q$-invariant 
if and only if $N^-$ is orientable, and in this case, it is the lifting of $w \in \Omega_{eq}((\T \times \R^{n-1})/ \Z_m)$ which 
generate the kernel of $L_s$. The other non-zero eigenvalues of $L_s$, being non-zero eigenvalues of $\tilde L_s$, go to infinity 
when $s$ and $\kappa$ go to infinity. This argument has the following conclusion 
\begin{lem}\label{consumo}
The relation
\[\lim_{s\to\infty}\sum_{o}\tr\,\phi_o\phi(\bar{L}_s^k)=d_k\]
holds, where $o$ runs over the critical orbits of the Morse function $f$ and $d_k$ 
denotes the number of critical orbits with index $k$ having orientable bundle $N^-$ 
(and $d_k=0$ for $k\geq n$).
\end{lem} 

\noindent Now we have everything to prove the main theorem \ref{thm1}. \\
\begin{pf}[{\bf of the main theorem \ref{thm1}}]
The claimed inequalities follow  
directly from theorem \ref{thm2} and lemmas \ref{keylemma}, \ref{consump} and \ref{consumo}. 
We just need to show that for $k\geq n+1$ the Morse inequalities do not provide new ones  
and reduce to lower order Morse inequalities. We do this for $k=n$, the general 
case is similar and follows from \eqref{2eqper}. 
For this purpose note that $d_k=0$ for $k\geq n$, therefore, by 
\ref{keylemma} and \ref{consump} we have 

\[\lim_{s\to\infty}(\mu_s^{n+1}-\mu_s^n)=(c_{n-1}+c_{n-3}+\dots)-(c_n+c_{n-2}+\dots)\]
The right side of this equality is $(-1)^{n-1}$ times the sum of the indices of the vector 
field $v$ on its singularities which equals $(-1)^{n-1}\chi(M)$, by the Poincare-Hopf theorem. 
This and lemma \ref{localEu} show that 
$\tilde c_{n+1}-\tilde c_n=\beta_{eq}^{n+1}-\beta_{eq}^{n}$ and complete the proof.
\end{pf}

\bibliographystyle{plain}

\end{document}